\documentstyle[11pt, pb-diagram]{article}
\input amssym.def
\newenvironment{Literatur}[1]{%
\begin{thebibliography}{#1}}{\end{thebibliography}}%

\title{Some Remarks on the Hijazi Inequality and Generalizations of the Killing Equation for Spinors.  \footnote{Supported by the SFB 288 and the Graduiertenkolleg ''Geometrie und Nichtlineare Analysis'' of the DFG.}}
\date{}

\begin{document}

%==========================
\newcommand{\D}{\displaystyle}
\newcommand{\upsp}{\phantom{l}}
\newcommand{\downsp}{\phantom{q}}
%==========================

\begin{center}
\maketitle

\vspace{-1.5cm}

{\large Thomas Friedrich\footnote[1]{E-mail: friedric@mathematik.hu-berlin.de} and Eui Chul Kim\footnote[2]{E-mail: zgast8@mathematik.hu-berlin.de} }\\
\mbox{}\\
{\small \it Humboldt-Universit\"at zu Berlin, Institut f\"ur Reine Mathematik,\\
Ziegelstra\ss e 13a, D-10099 Berlin, Germany}\\
\mbox{}\\
\today \\%Received: \\
\end{center}

\mbox{} \hrulefill \mbox{}\\

\newcommand{\vol}{\mbox{vol} \, }
\newcommand{\grad}{\mbox{grad} \, }

\begin{abstract} We generalize the well-known lower  estimates for the first eigenvalue of the Dirac operator on a compact Riemannian spin manifold proved by Th. Friedrich (1980) and O. Hijazi (1986, 1992). The special solutions of the Einstein-Dirac equation constructed recently by Friedrich/Kim are examples for the limiting  case of these inequalities. The discussion of the limiting case of these estimates yields two new field equations generalizing the Killing equation as well as the weak Killing equation for spinor fields. Finally, we discuss the 2- and 3-dimensional case in more detail. 
\end{abstract}

{\small
{\it Subj. Class.:} Differential Geometry.
 }\\

\setcounter{section}{0}

\mbox{} \hrulefill \mbox{}\\

\section{Introduction}

The first sharp estimate for the eigenvalues $\lambda$ of the Dirac operator defined on a compact $n$-dimensional Riemannian spin manifold was obtained by Th. Friedrich in 1980. Using a suitable deformation of the Riemannian connection he proved the inequality
$$ \lambda^2 \ge \frac{1}{4} \, \frac{n}{n-1} S_o , $$

where $S_o$ denotes the minimum of the scalar curvature $S$. The discussion of the limiting case yields a  stronger first order equation, the so-called spinorial Killing equation (see [5]). In 1986 O. Hijazi (see [9]) generalized this inequality. He combined the technique used before with a conformal change of the metric and obtained the inequality 
$$ \lambda^2 \ge \frac{1}{4} \, \frac{n}{n-1} \inf_M \Big\{ S + \frac{4(n-1)}{n-2} u^{-1} \Delta (u) \Big\} . $$

This estimate holds for any positive function $u$ defined on a compact Riemannian spin manifold of dimension $n \ge 3$ and the discussion of the limiting case yields again essentially the spinorial Killing equation. In particular, $\lambda^2$ is bounded from below by $\frac{1}{4} \, \frac{n}{n-1} \mu_1$, where $\mu_1$ is the first eigenvalue of the Yamabe operator. The simple identity
$$ \frac{n-1}{n-2} u^{-1} \Delta (u) = \Delta (f) - \frac{n-2}{n-1} |df|^2 \quad , \quad f= \frac{n-1}{n-2} \log (u) $$

allows us to rewrite the Hijazi inequality in the following equivalent form:
$$ \lambda^2 \ge \frac{n}{4(n-1)} \inf_M \Big\{ S + 4 \Delta (f) - 4 \frac{n-2}{n-1} |df|^2 \Big\} . $$

The advantage of this reformulation of the Hijazi inequality is the fact that the latter estimate is even true in dimension $n \ge 2$. This observation was made by C. B\"ar in 1992 (see [2]). In particular, using the Gauss-Bonnet Theorem one obtains for any metric $g$ on the sphere $S^2$ the inequality 
$$ \lambda^2  \ge \frac{4 \pi}{\vol (S^2,g)} . $$

Consequently, the Lott constant (see [14]) of the 2-sphere equals $4 \pi$, a value already conjectured by J. Lott in 1986. In 1992 O. Hijazi (see [10]) generalized the mentioned inequality by taking into account the energy-momentum tensor of the eigenspinor. This tensor occurs in the Einstein-Dirac equations describing  the interaction of a particle of spin $1/2$ with the gravitational field. Recently, we have constructed many solutions of this non-linear system (see [11]). Therefore, we have revisited the Hijazi inequality once again. Using deformations of the Riemannian connection depending on more free parameters then in all investigations before, we will prove an estimate (see Theorem 3.4) containing all these inequalities as special cases. Moreover, the weak Killing spinors that are special solutions of the Einstein-Dirac equation constructed in [11] are examples of spinors realizing  the limiting case in our new estimate. The discussion of the limiting case of these estimates yields two new field equations generalizing the spinorial  Killing equation as well as the weak Killing equation. \\

Both authors thank Ilka Agricola for her helpful comments and Heike Pahlisch for her competent and efficient \LaTeX work. We thank O. Hijazi for reading  a preliminary version of this paper and pointing out a mistake. During this discussion B. Ammann and C. B\"ar constructed examples showing that in contrast to the case $n=2$ in dimensions $n \ge 3$ the eigenvalue $\lambda^2$ is not bounded by the mean value of the scalar curvature.\\

\section{Generalization of the Hijazi inequality}

\vskip 0.3 cm 
Let $(M^n , g)$ be an $n$-dimensional connected oriented compact Riemannian spin manifold without boundary. We consider two conformally related metrics  $\overline{g} =  e^{- 2 h} g$,  where $h$ is a real-valued smooth function on $M^n$. Let us denote by $\Sigma (M)_g$ and $\Sigma (M)_{\overline{g}}$ the spinor bundle of $(M^n , g)$ and $(M^n , \overline{g})$,
respectively. There are natural isomorphisms  $ j : T(M) \longrightarrow
T(M)$ and  ${j} : \Sigma (M)_g \longrightarrow \Sigma (M)_{\overline{g}} $  preserving the inner products of vectors and spinors as well as the Clifford multiplication:
\begin{eqnarray*}
&    & \overline{g} (j X   ,   j Y) =  g(X ,  Y)  \quad , \quad < {j}   \varphi   ,   {j}   \psi >_{\overline{g}} \, \, = \, \, < \varphi ,   \psi >_g ,     \cr
&    &     \cr
&   & (j X) \cdot ({j}   \psi) = {j}   (X \cdot \psi) ,
\qquad X , Y  \in \Gamma (T(M)) , \, \, \varphi , \psi \in  \Gamma (\Sigma {(M)_g}).
\end{eqnarray*}

Let $\psi$ be a spinor field on $(M^n , g)$ and denote by $\overline{\psi} : = {j} (\psi)$ the corresponding spinor field on $(M^n, \overline{g})$. We use the same notation for vector fields, $\overline{X} := j(X)$. Then the following formulas relating the connections $\nabla, \overline{\nabla}$ as well as the Dirac operators $D, \overline{D}$ are well-known (see [3]):\\

\vskip 0.3 cm 
{\bf Lemma 2.1.}   
\begin{itemize}
\item[(i)] $\displaystyle \overline{{\rm grad}   \, }(f)   =   e^h   \overline{{\rm grad}  \, (f)} ,  $\\

\item[(ii)] $\displaystyle {\overline{\nabla}}_X (e^{{{n-1} \over 2} h}   \overline{\psi}) = \
 e^{{{n-1} \over 2} h}   \{   \overline{ \nabla_X \psi }   +   {n \over 2}   \overline{g} (\overline{{\rm grad}  \, }(h) , X) \overline{\psi}    +   {1 \over 2}   X \cdot  \overline{{\rm grad} \, }(h) \cdot
 \overline{\psi}   \}   ,   $\\

\item[(iii)] $\displaystyle \overline{D} (e^{{{n-1} \over 2} h}   \overline{\psi})  = \
e^{{{n+1} \over 2} h}   \overline{D \psi} ,  $\\

\item[(iv)] $\displaystyle (\overline{D} \circ \overline{D}) (e^{{{n-1} \over 2} h}   \overline{\psi}) = \
e^{{{n+3} \over 2} h}   \{   \overline{ D^2 \psi}   +   \overline{ {\rm grad} \, (h) \cdot D \psi }   \} .$
\end{itemize} 

\vskip 0.6 cm 
For later reference, we shall first assemble a few technical results. First, Lemma 2.1 (ii) implies immediately the following result.\\

\vskip 0.3 cm 
{\bf Corollary 2.2.}   {\it Suppose that  a spinor field $\overline{\psi}$ on $(M^n , \overline{g})$ satisfies the equation
$$ {\overline{\nabla}}_X \Big(e^{{{n-1} \over 2} h}   \overline{\psi}\Big) = e^{{{n-1} \over 2} h}   \Big\{ - {\lambda \over n} e^h   X \cdot \overline{\psi}   +   a   \overline{g} (\overline{{\rm grad} \, }(h) , X)    \overline{\psi}    +   b   X \cdot  \overline{{\rm grad} \, }(h) \cdot  \overline{\psi}   \Big\} $$

for some real numbers $ \lambda , a , b \in {\Bbb R}$ and for all vector fields $X$. Then the corresponding spinor field $\psi$ on $(M^n , g)$ satisfies the equations
$$ \nabla_X \psi = - {\lambda \over n}   X \cdot \psi   + \Big(a - {n \over 2}\Big) dh(X) \psi   +  
\Big(b - {1 \over 2}\Big) X \cdot {\rm grad} \, (h) \cdot \psi  $$  

and
$$ D \psi = \lambda \psi   +   (a - nb) {\rm grad} \, (h) \cdot \psi . \hspace{4.2cm} \mbox{} $$   }

\vskip 0.3 cm 
{\bf Lemma 2.3.}   {\it Let $\psi$ be an eigenspinor of the Dirac operator $D$ of $(M^n , g)$ with eigenvalue $0 \not= \lambda \in {\Bbb R}$ such that
$$ \nabla_X \psi =  -  {\lambda \over n}   X \cdot \psi   +   n d f (X) \psi   +    X \cdot {\rm grad} \,  (f) \cdot \psi $$
holds for some real-valued function $f$ and for all vector fields $X$. Then $f$ is constant
on $M^n$.

\vskip 0.3 cm 
Proof.}   The formula  $ {1 \over 2}   \mbox{Ric} (X) \cdot \varphi = D(\nabla_X  \varphi)   -   \nabla_X (D \varphi)   -  \sum\limits_{u=1}^n {E_u \cdot \nabla_{\nabla_{E_u} X}   \varphi } $ (see [11])
yields
\begin{eqnarray*}
\mbox{Ric} (E_l) \cdot \psi & = & \Big\{   {{4(n-1) \lambda^2} \over n^2} - 2 \triangle (f) + 4(n-2) \vert df \vert^2   \Big\} E_l \cdot \psi   +   2(n-2) \nabla_{E_l} df \cdot \psi     \cr
&    &       \cr                                                                                                                                                                              &    &  -   4(n-2) df(E_l) df \cdot \psi   -  
{{4 \lambda} \over n} E_l \cdot df \cdot \psi   +   {{4(n-2) \lambda} \over n} df(E_l) \psi \, \, .
\end{eqnarray*}

Contracting this equation we obtain
$$ S \psi = \Big\{   {{4(n-1) \lambda^2} \over n} - 4(n-1) \triangle (f) + 4(n-1)(n-2) \vert df \vert^2   \Big\} \psi  -   {{8(n-1) \lambda} \over n}   df \cdot \psi  ,$$
and, consequently, $df \equiv 0$.         \hfill{Q.E.D.}\\

{\bf Lemma 2.4.} (see [11]) {\it Let $\psi$ be a spinor field on $(M^n , g)$ such that
\[ \nabla_X \psi   =   \beta(X) \cdot \psi   +   n df(X) \psi   +   X \cdot {\rm grad} \, (f) \cdot \psi \]

\medskip

holds for a real-valued function $f$, for a symmetric (1,1)-tensor field $\beta$ and for all vector fields $X$. Then we have the formula}
\[ \{ {\rm Tr}(\beta) \}^2   =   { S \over 4} +  \vert \beta \vert^2 + (n-1) \triangle (f) - (n-1)(n-2) \vert df \vert^2  . \]

\medskip

The scalar curvature $\overline{S}$ of  $(M^n , \overline{g})$ is related to the
scalar curvature $S$ of $(M^n , g)$ by a well-known formula. We formulate this equation in two different ways, first for all dimensions $n \ge 2$ and then for dimensions $n \ge 3$:\\

\vskip 0.3 cm 
{\bf Lemma 2.5.}  
\begin{eqnarray*}
\overline{S} & = & e^{2h}   \Big\{   S   -   2(n-1) \triangle (h)   -   (n-1)(n-2) \vert dh \vert^2   \Big\} \quad                                     (\mbox{if} \, \, n \geq 2)    \cr
&     &     \cr
& = & e^{2h}   \Big\{   S   +   { {4(n-1)} \over {n-2}}   e^{{{n-2} \over 2} h}   \triangle
(e^{- {{n-2} \over 2} h})   \Big\}  \hspace{2.1cm} (\mbox{if} \,\,  n \geq 3) .
\end{eqnarray*}

\vskip 0.3 cm 
Let us now proceed to the main topic of this article. For this, we repeat once again the proof of the Hijazi inequality.  Let us denote  the real part of the hermitian inner product of spinors  by $(,)_g$ (resp. $(,)_{\overline{g}}$) and  let $v_g$ (resp. $v_{\overline{g}}$) be the volume form of  $(M^n , g)$ (resp. $(M^n , \overline{g})$). For shortness we introduce the notation $\psi_h :=  e^{{{n-1} \over 2} h}   \psi$. Using Lemma 2.1 and 2.5 one verifies:
\begin{eqnarray*}
0 & \le & \sum_{i=1}^n    \int_{M^n \downsp}   \Big(\overline{\nabla}_{\overline{E}_i} \overline{\psi}_h + {\lambda \over n} e^h   \overline{E}_i \cdot  \overline{\psi}_h   ,   \overline{\nabla}_{\overline{E}_i} \overline{\psi}_h + {\lambda \over n} e^h   \overline{E}_i \cdot  \overline{\psi}_h  \Big)_{\overline{g}}  \,  v_{\overline{g}}          \cr
&    &     \cr
& = & \int_{M^n \downsp}  \Big\{   \Big((\overline{D} \circ \overline{D}) (\overline{\psi}_h) - {1 \over 4}                                     \overline{S}  \,  \overline{\psi}_h   ,  
 \overline{\psi}_h\Big)_{\overline{g}}   -   {{2 \lambda} \over n} e^h  \Big(\overline{D} (\overline{\psi}_h  ) , \overline{\psi}_h \Big)_{\overline{g}}   +   {{\lambda^2} \over n} e^{2h} \Big(\overline{\psi}_h ,  \overline{\psi}_h \Big)_{\overline{g}}     \Big\} v_{\overline{g}}     \cr
&    &     \cr
& = &  \int_{M^n \downsp} e^{(n+1)h} \Big(\lambda^2 - {1 \over 4} e^{-2h}    \overline{S} - {{2 \lambda^2} \over n}+ {\lambda^2 \over n}  \Big) \Big(\psi , \psi \Big)_g   e^{-nh} v_g       \cr
&    &     \cr
& = & \int_{M^n \downsp} e^h \Big\{   {{n-1} \over n} \lambda^2  - {S \over 4}  + {{n-1} \over 2} \triangle(h) + {{(n-1)(n-2)} \over 4} \vert dh \vert^2   \Big\} \Big(\psi , \psi \Big)_g v_g     .     \cr
&    &     \cr
\end{eqnarray*}

If the dimension satisfies $n \ge 3$, we rewrite the latter equation in the following equivalent form:
\begin{eqnarray*}
0 & \le & \int_{M^n \downsp} e^h \Big\{   {{n-1} \over n} \lambda^2  - {S \over 4}  - {{n-1} \over {n-2}}
u^{-1}  \triangle (u)   \Big\} \Big(\psi , \psi \Big)_g \, v_g      , \cr
\end{eqnarray*}

where the arbitrary positive function $u$ is related to $h$ by  $u : =  e^{- {{n-2} \over 2} h}$. Then we obtain the Hijazi inequality (see [9]) 
$$ (\ast) \qquad  \lambda^2 \geq {n \over {4(n-1)}} \inf_M \Big\{   S  +  {{4(n-1)} \over {n-2}} u^{-1} \triangle(u)   \Big\} \quad   , \quad n \ge 3 . $$

\vskip 0.3 cm 
The fourth line of the above calculation yields an equivalent version of this  inequality  valid for all dimensions $n \ge 2$
\[  \lambda^2 \geq {n \over {n-1}} \inf_M \Big\{   {S \over 4}  +  \triangle(f) - {{n-2} \over {n-1}} \vert df \vert^2   \Big\} . \]

An eigenspinor $\psi$ of the Dirac operator $D$ for the limiting eigenvalue $\lambda$ satisfies,  by Corollary 2.2,  the stronger field equation  
\[ \nabla_X \psi = - {\lambda \over n}   X \cdot \psi   -   {n \over 2} dh(X) \psi  
 -   {1 \over 2} X \cdot {\rm grad} \, (h) \cdot \psi   \]

\medskip
  
for all vector fields $X$. In case $\lambda \neq 0$, we conclude, by Lemma 2.3,  that $h$ is constant and,  therefore,  $\psi$ is a Killing spinor. In case $\lambda = 0$, the spinor   $e^{{{n-1} \over 2} h}   \overline{\psi}$ is parallel on $(M^n , \overline{g})$. \\

\vskip 0.3 cm 
Now we prove the main result of this section. \\

{\bf Theorem 2.6.}   {\it Let $(M^n , g)$ be a compact Riemannian spin manifold of dimension $n \geq 2$. For any eigenvalue $\lambda$ of the Dirac operator $D$ the inequality
\[ \lambda^2 \geq {n \over {n-1}} \inf_M \Big\{  {S \over 4}  - \Big({{n-1} \over 2}- a \Big) \triangle(h)  
- \Big( {{(n-1)(n-2)} \over 4} + a^2 + nb^2 - na + 2a - 2ab \Big) \vert dh \vert^2 \Big\}   \]

\medskip

holds for all real-valued functions $h$ on $M^n$ and for all real numbers $a , b \in {\Bbb R}$. Equality occurs if and only if either $(M^n , g)$ admits a Killing spinor or if there is a conformally equivalent metric $\overline{g}$ such that $(M^n , \overline{g})$ admits  a parallel spinor. 

\vskip 0.3 cm 
Proof.}   Let $\psi$ be an eigenspinor of the Dirac operator $D$ with eigenvalue $\lambda$ and consider again the spinor field $\psi_h =  e^{{{n-1} \over 2} h}   \psi$. Define for any real numbers $a$ and $b$ the spinor
$$ \overline{P} (\overline{X}) := \overline{\nabla}_{\overline{X}} \overline{\psi}_h + {\lambda \over n} e^h   \overline{X} \cdot  \overline{\psi}_h - a e^h dh(X)  \overline{\psi}_h - b e^h  \overline{X} \cdot \overline{{\rm grad} \, (h)} \cdot \overline{\psi}_h  . $$

Then a direct calculation using Lemma 2.1, Lemma 2.5 and the Schr\"odinger-Lichnerowicz formula $D^2 = \Delta + \frac{1}{4} S$ yields the equation 
$$ 0  \le  \sum\limits^{n}_{i=1} \int_{M^n} (\overline{P} (\overline{E}_i) \, , \, \overline{P} (\overline{E}_i)) {v_{\overline{g}}}  = \int_{M^n} e^h H | \psi |^2 v_g  , $$

where the function $H$ is given by the formula
$$ H=     {{n-1} \over n} \lambda^2 - { S \over 4}  + \Big({{n-1} \over 2}- a \Big) \triangle(h)  +   \Big({{(n-1)(n-2)} \over 4} + a^2 + nb^2 - na + 2a - 2ab \Big) \vert dh \vert^2   . $$

This identity proves the inequality of the theorem. Now we discuss the limiting case. By Corollary 2.2  we obtain the following differential equation for an eigenspinor $\psi_1$ with the limiting eigenvalue $\lambda_1$
\[ \nabla_X \psi_1 = - {\lambda_1 \over n}   X \cdot \psi_1   + \Big(a - {n \over 2}\Big) dh_1(X) \psi_1   +   \Big(b - {1 \over 2}\Big) X \cdot {\rm grad} \, (h_1) \cdot \psi_1  ,  \]

\medskip

as well as the condition  that  $(a-nb) {\rm grad} \, (h_1) \equiv 0$.   If $\lambda_1 \neq 0$ and  $h_1$ is not constant, then $a = nb$ and we conclude $b= \frac{1}{2}$ by Lemma 2.3, i.e., $\psi_1$ is a Killing spinor. In case that $\lambda_1 = 0$ and $h_1$ is not constant, we have
\[ \nabla_X \psi_1 =  n \Big(b - {1 \over 2}\Big) dh_1(X) \psi_1   +   
\Big(b - {1 \over 2}\Big) X \cdot {\rm grad} \, (h_1) \cdot \psi_1  . \]

\medskip

Moreover,  the limiting case of the inequality as well as Lemma 2.4 yield the two equations
\begin{eqnarray*}
0 & = &  {S \over 4}  - \Big({{n-1} \over 2}- nb\Big) \triangle(h_1) 
- \Big( {{(n-1)(n-2)} \over 4} + n(n-1)b^2- n(n-2)b \Big) \vert dh_1 \vert^2     \cr
&    &     \cr
& = & {S \over 4} + (n-1)\Big(b - {1 \over 2}\Big) \triangle(h_1) - (n-1)(n-2)\Big( b - {1 \over 2}\Big)^2 \vert dh_1 \vert^2 \ .
\end{eqnarray*}

Therefore, $b=0$, thus implying by Corollary 2.2 that $\displaystyle {\overline{\nabla}}_X (e^{{{n-1} \over 2} h_1}  \overline{\psi_1}) = 0$. Consequently, $\displaystyle e^{{{n-1} \over 2} h_1} \, \overline{\psi_1}$ is a parallel spinor with respect to the metric $\overline{g} = e^{- 2 h_1} g$. The converse can easily be  proved using Lemma 2.4.          \hfill{Q.E.D.}\\

{\bf Remark.} One can maximize this estimate with respect to the quadratic term $|dh|^2$ only. Then the optimal parameters are
$$a = \frac{n}{2} \cdot \frac{n-2}{n-1} \quad , \quad b = \frac{1}{2} \cdot \frac{n-2}{n-1}$$

and we obtain
$$ \lambda^2 \ge \frac{n}{n-1} \inf_M \Big\{ \frac{S}{4} + \frac{1}{2(n-1)} \Delta (h) - \frac{n-2}{4(n-1)} |dh|^2 \Big\} . $$

\section{\mbox{Lower eigenvalue estimates using the energy-}momentum tensor}

\vskip 0.3cm 
Any eigenspinor $\psi$ of the Dirac operator $D$ of $(M^n , g)$ induces a symmetric (0,2)-tensor field $T_{\psi}$ defined by
$$ T_{\psi} (X , Y)   =   (X \cdot \nabla_Y \psi + Y \cdot \nabla_X \psi   ,   \psi  ) ,$$

which is the energy-momentum tensor in the Einstein-Dirac equation (see [11]). 
Over the open dense subset $M_{\psi} : = \{ x \in M^n : \psi (x) \neq 0 \}$ we define the tensor field 
$$ \widehat{T}_{\psi} (X , Y)   : =   (X \cdot \nabla_Y  \widehat{\psi} + Y \cdot \nabla_X \widehat{\psi}, \widehat{\psi}) =   {1 \over {\vert \psi \vert^2}} T_{\psi} (X,Y)  ,$$

where   $\displaystyle \widehat{\psi} : =  {\psi \over {\vert \psi \vert}}$ is the normalized spinor. Hijazi (see [10]) proved the following eigenvalue estimates for the Dirac operator depending on the scalar curvature $S$, the first eigenvalue $\mu_1$ of the Yamabe operator and on the length of $\widehat{T}_{\psi}$:
$$ (\ast\ast)  \qquad  \lambda^2   \geq   {1 \over 4} \inf_{M_{\psi}} (S   +   \vert  \widehat{T}_{\psi} \vert^2  )  \qquad \mbox{and} \qquad  \lambda^2   \geq   {1 \over 4} \mu_1   +   {1 \over 4}  \inf_{M_{\psi}}  \vert  \widehat{T}_{\psi} \vert^2  .$$

In this section we will improve the inequalities $(\ast\ast)$ and show the limiting case explicitly. For $g, \overline{g}$ and $\psi_h$ defined as above, one easily verifies the following formulas:

\vskip 0.3 cm 
{\bf Lemma 3.1.}  
\begin{itemize}
\item[(i)] $\displaystyle {\widehat{T}}_{\overline{\psi}_h}(\overline{E}_k  , \overline{E}_l ) = e^h {\widehat{T}}_{\psi} (E_k , E_l)   \quad  (1 \leq k , l \leq n)  $  ,   \\
\item[(ii)] $\displaystyle  \sum_{i=1}^n   ({\widehat{T}}_{\overline{\psi}_h} (\overline{E}_i) \cdot \overline{\nabla}_{\overline{E}_i} \overline{\psi}_h   ,   \overline{\psi}_h  )_{\overline{g}} = {1 \over 2} e^{(n+1)h} \vert {\widehat{T}}_{\psi} \vert^2  (\psi , \psi)_g .$
\end{itemize}   

\vskip 0.3 cm 
Corollary 2.2 can easily be generalized:

\vskip 0.5 cm 
{\bf Lemma 3.2.} \,  {\it Let $\overline{\psi}$ be a spinor field on $(M^n , \overline{g})$ and $U$ an open subset of $M^n$.    Suppose that the spinor field $\overline{\psi}$ satisfies, on $U$, the equation
$$ {\overline{\nabla}}_X \Big( e^{{{n-1} \over 2} h}  \overline{\psi}\Big)  = 
 e^{{{n-1} \over 2} h}  \Big\{  e^h  \beta(X) \cdot \overline{\psi}  +  a  \overline{g} ( \overline{{\rm grad}}(h) , X )   \overline{\psi}   +  b X \cdot  \overline{{\rm grad}}(h) \cdot \overline{\psi}  \Big\}  $$

for some symmetric (1,1)-tensor field $\beta$, a real-valued function $h$ and for                                                            all real numbers  $ \lambda , a , b \in {\Bbb R}$ ( $\beta$ and $h$ may be defined on the subset $U$ only).  Then the corresponding spinor field $\psi$ on $(M^n , g )$ satisfies, on $U$, the equations
$$ \nabla_X \psi  =  \beta(X) \cdot \psi  + \Big(a - {n \over 2}\Big) dh(X) \psi  +  
\Big( b - {1 \over 2}\Big) X \cdot {\rm grad}(h) \cdot \psi  $$   
and 
$$ D \psi  =  - {\rm Tr}(\beta) \psi  +  ( a - nb ) {\rm grad}(h) \cdot \psi .$$   }

\vskip 0.3 cm 
{\bf Lemma 3.3.} (see [11])    {\it Let $\psi$ be a non-trivial spinor field on $(M^n , g)$ such that, on a connected  open subset $U \subset M^n$, the equation 
$$ \nabla_X \psi = \beta(X) \cdot \psi   +   n \alpha(X) \psi   +   X \cdot \alpha \cdot \psi $$

holds for a 1-form $\alpha$, a symmetric $(1,1)$-tensor field $\beta$ and for all vector fields $X$.   Then $\psi$ has no zeros in $U$ and $\alpha$ as well as $\beta$ are uniquely determined by the spinor field $\psi$ via the relations:
$$\alpha = {{d(\vert \psi \vert^2)} \over {2(n-1) \vert \psi \vert^2}} \qquad and \qquad \beta = -  \frac{1}{2} {\widehat{T}_{\psi}} . $$

In particular, the 1-form $\alpha$ is exact. }\\

{\bf Theorem 3.4.} \, {\it  Let $( M^n , g )$ be a compact Riemannian spin manifold of dimension $n \geq 2$. For any eigenspinor $\psi$ of the Dirac operator $D$ with eigenvalue $\lambda$ the inequality
$$  \lambda^2  \geq  \inf_{M_{\psi}} \Big\{  {S \over 4} + {1 \over 4} \vert {\widehat{T}_{\psi}} \vert^2  - \Big({{n-1} \over 2}- a\Big) \triangle(h)  - \Big(  {{(n-1)(n-2)} \over 4} + a^2 + nb^2 - na + 2a - 2ab \Big) \vert dh \vert^2  \Big\} $$

holds for all real-valued functions $h$ on $M^n$ and for all real numbers $a , b \in {\Bbb R}$. Equality occurs if and only if there exists an eigenspinor $\psi_1$ without zeros such that the equation
$$ \nabla_X \psi_1  =  - {1 \over 2} \widehat{T}_{\psi_1}(X) \cdot \psi_1  -  {n \over 2} dh_1(X) \psi_1  -  {1 \over 2} X \cdot {\rm grad}(h_1) \cdot \psi_1 $$

holds for all vector fields $X$ on $M^n$. In this limiting case the function $h_1$ is uniquely determined up to constants by the spinor field $\psi_1$ via the relation 
$$ h_1  =  -  {{ \log(\vert \psi_1 \vert^2) } \over {n-1}} .$$

Proof.}  We use a slight modification of the field $\overline{P} (\overline{X})$ used in the proof of Theorem 2.6. Namely, set 
\begin{eqnarray*}
\overline{Q} (\overline{X}) &  := & \overline{\nabla}_{\overline{X}} \overline{\psi}_h + \frac{1}{2} \widehat{T}_{\overline{\psi}_h} (\overline{X})  \cdot  \overline{\psi}_h - a e^h dh(X)  \overline{\psi}_h - b e^h  \overline{X} \cdot \overline{{\rm grad} \, (h)} \cdot \overline{\psi}_h   .
\end{eqnarray*}

Then one shows, as before, the equation 
$$ 0  \le  \sum\limits^{n}_{i=1} \int_{M^n} (\overline{Q} (\overline{E}_i)  ,  \overline{Q} (\overline{E}_i)) {v_{\overline{g}}} = \int_{M^n} e^h H | \psi |^2 v_g ,$$

where the function $H$ is now given by the formula
\[ H= \lambda^2 - { S \over 4} - {1 \over 4} \vert {\widehat{T}_{\psi}} \vert^2 + \Big({{n-1} \over 2}- a\Big) \triangle(h)   +   \Big({{(n-1)(n-2)} \over 4} + a^2 + nb^2 - na + 2a - 2ab \Big) \vert dh \vert^2  .\]

This yields the inequality of the theorem. For the limiting case,  we obtain by Lemma 3.2 the following differential equation for an eigenspinor $\psi_1$
$$ \nabla_X \psi_1  =   - {1 \over 2} \widehat{T}_{\psi_1}(X) \cdot \psi_1  + \Big(a - {n \over 2}\Big) dh_1(X) \psi_1  + \Big( b - {1 \over 2}\Big) X \cdot {\rm grad}(h_1) \cdot \psi_1   ,$$

as well as the condition $(a-nb) {\rm grad}(h_1) \equiv 0$ for a real-valued function $h_1$ defined on the whole manifold $M^n$. According to Lemma 3.3 the eigenspinor $\psi_1$ does not vanish anywhere. The trace of  $\widehat{T}_{\psi_1}$ is related to the eigenvalue by ${\rm Tr}(\widehat{T}_{\psi_1}) = 2 \lambda_1$.  In case that $h_1$ is not constant, we have $a = nb$ and, consequently,
$$\nabla_X \psi_1  =  - {1 \over 2} \widehat{T}_{\psi_1}(X)  \cdot \psi_1  + n \Big( b - {1 \over 2}\Big) dh_1(X) \psi_1  + \Big( b - {1 \over 2}\Big) X \cdot {\rm grad}(h_1) \cdot \psi_1  .$$

From the limiting case of the inequality of the theorem and Lemma 2.4 we obtain the equation
\begin{eqnarray*}
\lambda_1^2 & = &  {S \over 4} + {1 \over 4} \vert {\widehat{T}_{\psi_1}} \vert^2 - \Big({{n-1} \over 2}- nb\Big) \triangle(h_1)  - \Big\{ {{(n-1)(n-2)} \over 4} + n(n-1)b^2- n(n-2)b \Big\} \vert dh_1 \vert^2     \cr
&    &     \cr
& = & {S \over 4} + \frac{1}{4} \vert \widehat{T}_{\psi_1} \vert^2 + (n-1)\Big(b - {1 \over 2}\Big) \triangle(h_1) - (n-1)(n-2)\Big(b - {1 \over 2}\Big)^2 \vert dh_1 \vert^2 ,
\end{eqnarray*}

i.e., $b = 0$. \hfill{Q.E.D.}\\

\vspace{0.3cm}

{\bf Corollary 3.5.}   {\it  For any positive function $u$ the inequality
$$  \lambda^2 \geq {1 \over 4} \inf_{M_{\psi}} \Big\{   S  +  \vert {\widehat{T}_{\psi}} \vert^2 + {{4(n-1)} \over {n-2}} u^{-1} \triangle(u)  \Big\}  $$

holds. In particular, if $u$ is the eigenfunction for the first eigenvalue $\mu_1$ of the Yamabe operator $ L = {{4(n-1)} \over {n-2}} \triangle + S$ we obtain the inequality $(\ast \ast)$:
$$  \lambda^2 \geq {1 \over 4}  \mu_1   +   {1 \over 4} \inf_{M_{\psi}}  \vert {\widehat{T}_{\psi}} \vert^2 .$$ }

\vspace{0.3cm}

{\bf Remark.} One can again maximize the estimate of Theorem 3.4 with respect to the quadratic term $|dh|^2$. The optimal parameters are
$$a = \frac{n}{2} \cdot \frac{n-2}{n-1} \quad , \quad b = \frac{1}{2} \cdot \frac{n-2}{n-1}$$

and we obtain
$$ \lambda^2 \ge  \inf_{M_{\psi}} \Big\{ \frac{S}{4} + \frac{1}{4} |\widehat{T}_{\psi}|^2 + \frac{1}{2(n-1)} \Delta (h) - \frac{n-2}{4(n-1)} |dh|^2 \Big\} . $$

\vspace{0.3cm}

{\bf Remark.} Because of  the relation 

$$   g \Big({\widehat{T}}_{\psi}  -  {{2 \lambda} \over n} g   ,   {\widehat{T}}_{\psi}  -  {{2 \lambda} \over n} g \Big)   =  \vert {\widehat{T}}_{\psi} \vert^2  - {{4 \lambda^2} \over n}   \geq   0  $$

the inequality of Theorem 3.4 is stronger than all the estimates $(\ast), (\ast \ast)$ and the estimate in Theorem 2.6.\\

\vspace{0.3cm}

{\bf Example 1.}  Let $\psi_1$ be a weak Killing spinor on $( M^n , g )$ with WK-number $\lambda_1$, i.e., a solution of the following differential equation (see [11], $n \geq 3$)
$$ \nabla_X \psi  =  { n \over 2 \, (n-1)S }  dS(X)  \psi  +  { 2 \, \lambda \over (n-2)S } 
 \mbox{Ric} (X) \cdot \psi  -  { \lambda \over n-2 }  X  \cdot \psi  +  { 1 \over 2  (n-1)S } X \cdot dS \cdot \psi \ , $$ 

where $\lambda$ is a real number.  Then $\psi_1$ satisfies the limiting case of the inequality of Theorem 3.4. Indeed, we have
$$  \widehat{T}_{\psi_1}  =  - { 4 \lambda_1 \over (n-2) S }  \mbox{Ric}  +  {2 \lambda_1 \over n-2 } {\rm Id}  \qquad \mbox{and} \qquad h_1  = - {{\log (\vert S \vert)} \over {n-1}} .$$

{\bf Example 2.} (see [11])   Let $(M^{2m+1}, \phi , \xi , \eta , g)$ be a simply connected Sasakian spin manifold with Ricci tensor $\mbox{Ric} = (2m - 4b) g + 4b \eta \otimes \eta , \, \, b \in {\Bbb R}$. Then there exists a non-trivial eigenspinor $\psi_1$ of the Dirac operator with eigenvalue $\lambda_1 = m + {1 \over 2} - b$ such that
$$ \nabla_X \psi_1 = - {1 \over 2} X \cdot \psi_1   +   b \eta(X) \xi \cdot \psi_1  .$$

This eigenspinor $\psi_1$ is an example of the limiting case of the inequality. Moreover, the length of its energy-momentum tensor is given by
$$ | \widehat{T}_{\psi_1} |^2 = \frac{4 \lambda_1^2 + 8mb^2}{2m+1} \ge \frac{4 \lambda_1^2}{2m+1} . $$

The discussion of the limiting case in the inequalities of Theorem 3.4  yields two new equations generalizing the Killing equation (see [5]) as well as the weak Killing equation (see [11]) for spinor fields.\\

\vspace{0.3cm}

{\bf Definition.} Let $(M^n,g)$ be a Riemannian spin manifold. A spinor field $\psi$ without zeros will be called
\begin{itemize}
\item[(i)] a {\it $T$-Killing spinor} if the trace ${\rm Tr} (\widehat{T}_{\psi}) = \frac{1}{|\psi|^2} {\rm Tr} (T_{\psi})$ is constant and $\psi$ is a solution of the equation
$$ \nabla_X \psi = - \frac{1}{2} \widehat{T}_{\psi} (X) \cdot \psi \quad , \quad X \in T(M^n) . $$

\item[(ii)] a {\it weak $T$-Killing spinor} if the trace ${\rm Tr} (\widehat{T}_{\psi}) = \frac{1}{|\psi|^2} {\rm Tr} (T_{\psi})$ is constant and $\psi$ is a solution of the equation
$$ \nabla_X \psi = - \frac{1}{2} \widehat{T}_{\psi} (X) \cdot \psi - \frac{n}{2} dh (X) \psi - \frac{1}{2} X \cdot \mbox{grad} \, (h) \cdot \psi , $$

where $h$ is the function $\displaystyle h= - \frac{\log (|\psi |^2)}{n-1}$.
\end{itemize}

%\smallskip

The following table lists the different kinds of eigenspinors of the Dirac operator we introduced as well as the necessary geometric condition for the underlying space.

\vspace{-0.1cm}

{\small
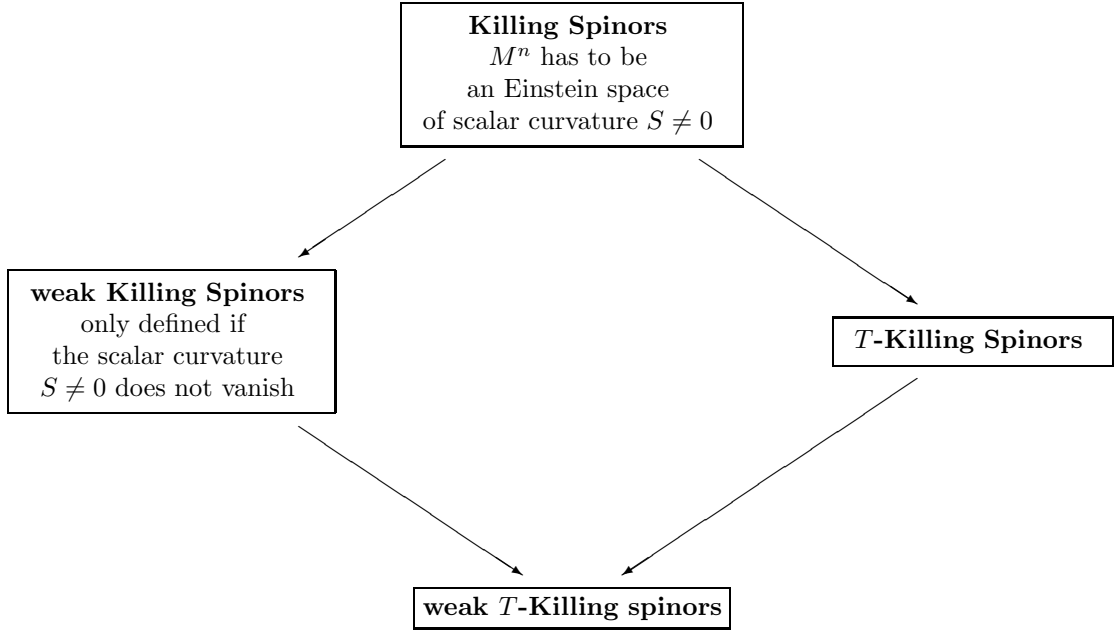
\[
\begin{diagram}
\node{} \node{\fbox{$\begin{array}{c}\mbox{\bf Killing Spinors}\\
					\mbox{$M^n$ has to be}\\					\mbox{an Einstein space}\\
					\mbox{of scalar curvature $S \not= 0$}
		\end{array}$
}
} \arrow{sw} \arrow{se} \node{}\node{}\\
\node{\fbox{$\begin{array}{c}\mbox{\bf weak Killing Spinors}\\
					\mbox{only defined if }\\
					\mbox{the scalar curvature}\\
					\mbox{$S \not= 0$ does not vanish}
		\end{array}$
}
} \arrow{se} \node{}  \node{\fbox{$\begin{array}{c}\mbox{\bf $T$-Killing Spinors}\\
		\end{array}$
}} \arrow{sw}\node{}\\
\node{}  \node{\fbox{\mbox{\bf weak $T$-Killing spinors}}} \node{}\node{}\\
\end{diagram}
\]}

\vspace{-2cm}

\section{The 2- and 3-dimensional case}

\vskip 0.3 cm \noindent
In this section we investigate the 2- and 3-dimensional case and present some properties of 
eigenspinors of the Dirac operator. For algebraic reasons  we can express, in these dimensions,  
the covariant derivative of an eigenspinor by the spinor itself  (see [11]):

\vskip 0.3 cm \noindent 
{\bf Lemma 4.1.}  {\it Let $( M^n , g )$ be a 2- or 3-dimensional Riemannian spin manifold,                                                             and let $\psi$ be an eigenspinor of the Dirac operator $D$ with eigenvalue $\lambda \in {\Bbb R}$. Then we have on the subset $M_{\psi}$
$$ \nabla_X \psi  =   - {1 \over 2} \widehat{T}_{\psi}(X) \cdot \psi +  n \alpha (X) \psi  +  X \cdot \alpha \cdot \psi $$

for a 1-form $\alpha$, which is uniquely determined by the spinor field $\psi$ via the relation 
$$  \alpha  =  {{d \{ \log(\vert \psi \vert^2) \}} \over {2(n-1)}} .$$  }

\vskip 0.3 cm \noindent
In any dimension we have proved (see [11]) the following estimate for the eigenvalue of the Dirac operator.

\vskip 0.3 cm \noindent
{\bf Lemma 4.2.}(see [11]) \, {\it Let $\psi$ be an eigenspinor of the Dirac operator $D$ with eigenvalue $\lambda \in {\Bbb R}$. Then the following inequality holds at any point of $M_{\psi}$ :
$$ \lambda^2 \quad \geq \quad {S \over 4}  +  {{\vert  \widehat{T}_{\psi} \vert^2} \over 4}  +  
 { {\triangle ( \vert \psi \vert^2 )} \over { 2 \, \vert \psi \vert^2} }  + {n \, {\vert d ( \vert \psi \vert^2 ) \vert^2} \over {4(n-1) \, \vert \psi \vert^4} } . $$  

Equality occurs if and only if there exists an eigenspinor $\psi_1$ of $D$  as well as a
1-form $\alpha_1$ such that on $M_{\psi_1}$
$$ \nabla_X \psi_1  =   - {1 \over 2} \widehat{T}_{\psi_1}(X) \cdot \psi_1 +  n \alpha_1 (X) \cdot \psi_1  +  X \cdot \alpha_1 \cdot \psi_1 $$ holds for all  vector fields $X$.   }

\vskip 0.3 cm \noindent
A direct consequence of Lemma 4.1 and 4.2 is the next

\vskip 0.3 cm \noindent
{\bf Theorem 4.3.}  {\it Let $( M^n , g )$ be a 2- or 3-dimensional Riemannian spin manifold                                                             and $\psi$ be an eigenspinor of the Dirac operator $D$ with eigenvalue $\lambda \in {\Bbb R}$. Then we have at any point of $M_{\psi}$
$$ \lambda^2 =  {S \over 4}  +  {{\vert  \widehat{T}_{\psi} \vert^2} \over 4}  +     
 { {\triangle ( \vert \psi \vert^2 )} \over { 2  \vert \psi \vert^2} } +  
{n  {\vert d ( \vert \psi \vert^2 ) \vert^2} \over {4(n-1)  \vert \psi \vert^4} } . $$

In particular, if both the scalar curvature $S$ and $\vert \psi \vert^2$ are constant, then ${\vert  \widehat{T}_{\psi} \vert^2}$ is constant.   }

\vskip 0.3 cm 
{\bf Theorem 4.4.}   {\it Let $(M^3 , g)$ be a 3-dimensional compact Riemannian spin manifold and $\psi$ a nowhere-vanishing eigenspinor of the Dirac operator $D$ with
eigenvalue $\lambda$. Then we have}
$$  \lambda^2 \cdot \mbox{vol} (M^3, g) \quad \leq \quad  {1 \over 4}   \int_{M^3}  \Big\{ S +  \vert  \widehat{T}_{\psi} \vert^2  \Big\} .$$

{\it Equality occurs if and only if $\vert \psi \vert^2$ is constant.

\vskip 0.4 cm 
Proof.}   By Theorem 4.3 we have
$$ \lambda^2   =   { S \over 4}   +   {{\vert  \widehat{T}_{\psi} \vert^2} \over 4}    +   2 \triangle (f)   -   2 \vert df \vert^2 , $$ where $\displaystyle  f = {1 \over 4} {{\rm log}(\vert \psi \vert^2)}$. Integrating this equation  we immediately obtain the inequality of the theorem.    \hfill{Q.E.D.}\\

{\bf Remark.} Let $(M^3,g)$ be a 3-dimensional Riemannian manifold with nowhere-vanishing scalar  curvature $S$. Let $\psi$ be an Einstein spinor (see [11]) for the eigenvalue $\lambda$, i.e.,  a solution of the non-linear system $D(\psi)= \lambda \psi$, $\mbox{Ric} - \frac{1}{2} Sg =  \pm \frac{1}{4} T_{\psi}$. Then
$$ | \widehat{T}_{\psi} |^2 = 4 \lambda^2 \Big( \frac{4 | \mbox{Ric} |^2}{S^2} - 1 \Big) $$

and, hence, Theorem 4.4 yields
$$ \lambda^2 \Big\{ \mbox{vol} (M^3,g) - 2 \int_{M^3} \frac{| \mbox{Ric}|^2}{S^2} \Big\} \le \frac{1}{8} \int_{M^3} S . $$

\vskip 0.3 cm 
{\bf Theorem 4.5.} {\it Let $(M^2 , g)$ be a 2-dimensional compact Riemannian spin manifold                                                             and $\psi$ a nowhere-vanishing eigenspinor of the Dirac operator $D$ with eigenvalue $\lambda$. Then we have }

\begin{itemize}
\item[(i)] $\displaystyle \lambda^2 =  {{\pi \chi({M}^2)} \over {{\rm vol}({M}^2 , g)}} 
+ {1 \over 4 {\rm vol} (M^2 , g)} \int_{M} \vert  \widehat{T}_{\psi} \vert^2 , $\\
\item[(ii)] $\displaystyle \int_M {\rm det} ( \widehat{T}_{\psi})  =  2 \pi \chi({M}^2)$ .
\end{itemize}

{\it Assume that $\vert \psi \vert^2$ is constant. Then $\widehat{T}_{\psi}$ is non-degenerate at a point 
$x \in M^2$ if and only if the scalar curvature $S$ does not vanish at $x$. 

\vskip 0.3 cm 
Proof.}    By Theorem 4.3 we have
$$ \lambda^2 = { S \over 4}   +   {{\vert  \widehat{T}_{\psi} \vert^2} \over 4}    +    \triangle (f)  , $$ 

where $\displaystyle  f = {1 \over 2} {{\rm log}(\vert \psi \vert^2)}$. The Gauss-Bonnet Theorem yields immediately the first equality of the theorem. Inserting

\begin{displaymath} {\rm det} ( \widehat{T}_{\psi}) =  {1 \over 2} \{ {\rm Tr}(\widehat{T}_{\psi}) \}^2-\frac{1}{2} {\rm Tr} \{ (\widehat{T}_{\psi})^2 \} = 2 \lambda^2 - \frac{1}{2} \vert \widehat{T}_{\psi} \vert^2 \end{displaymath}

into $\lambda^2 = { S \over 4} + {{\vert  \widehat{T}_{\psi} \vert^2} \over 4} + \triangle (f)$, we obtain
${\rm det} (\widehat{T}_{\psi}) = {S \over 2} + 2 \triangle (f)$, which implies the second identity (ii) as well as the last statement of the theorem.    \hfill{Q.E.D.}
\\

\vskip 0.3 cm
 
{\bf Remark. }  Since $\vert  \widehat{T}_{\psi} \vert^2 \geq 2 \lambda^2$ in the 2-dimensional case, Theorem 4.5 gives the ine\-quality
$$ \lambda^2 \geq  {{ 2 \pi \chi(M^2)} \over {{\rm vol}(M^2 , g)}} .$$

{\bf Example 1.} Let $f :M^2 \hookrightarrow {\Bbb R}^3$ be an isometric immersion of a closed surface $M^2$ into the Euclidean space ${\Bbb R}^3$ and suppose that the mean curvature $H$ is constant. A fixed parallel spinor $\Phi$ on ${\Bbb R}^3$ and its restriction onto the surface $M^2$ define an eigenspinor $\varphi$ of length one of the Dirac operator $D$ on the surface $(M^2,g)$. Moreover, this eigenspinor is a solution of the twistor type equation
$$ \nabla_X \varphi = - \frac{1}{2} \mbox{II} (X) \cdot \varphi \quad , \quad X \in T(M^2) , $$

where II denotes the second fundamental form of the surface (see [8]). The length $|T_{\varphi}|^2$ coincides with the length $|\mbox{II}|^2$ of the second fundamental form. The formulas of Theorem 4.5 are then simply the Gauss-Bonnet Theorem and
$$ H^2 = \frac{\pi \chi (M^2)}{\vol (M^2,g)} + \frac{1}{4 \vol (M^2,g)} \int_{M^2} |\mbox{II}|^2 . $$

Notice that this yields examples of $T$-Killing spinors on any surface of constant mean curvature in ${\Bbb R}^3$.\\

{\bf Example 2.} Consider the 2-dimensional torus $T^2 = {\Bbb R}^2 / {\Bbb Z}^2$ equipped with an $S^1$-invariant Riemannian metric
$$ g=h^4 (x) \{ dx^2 + dy^2 \}$$

and denote by $\lambda^2_1 (l)$ the first eigenvalue of the Dirac operator with respect to the trivial spin structure such that the $S^1$-representation of weight $l \not= 0$ occurs in the eigenspace $E (\lambda)$. Then the multiplicity of this representation in $E(\lambda)$ is one and the corresponding  unique eigenspinor  does not vanish at all (see [1]). \\

\vskip 1.0 cm  

\section{References}  

\vspace{-0.5cm}
\begin{Literatur}{MMM}
\bibitem[1]{} I. Agricola, B. Ammann and Th. Friedrich, A comparison of the eigenvalues of the Dirac and the Laplace operator on the two-dimensional torus, to appear in ''Manuscripta Mathematica''.
\bibitem[2]{} C. B\"{a}r, Lower eigenvalue estimates for Dirac operators, Math. Ann. 293 (1992), 39-46.   
\bibitem[3]{} H. Baum, Spin-Strukturen und Dirac-Operatoren \"{u}ber pseudo-Riemannschen Mannigfaltigkeiten, Teubner-Verlag, Leipzig 1981.
\bibitem[4]{} H. Baum, Th. Friedrich, R. Grunewald and I. Kath, Twistors and Killing spinors on Riemannian manifolds, Teubner-Verlag, Leipzig/Stuttgart 1991.   
\bibitem[5]{} Th. Friedrich, Der erste Eigenwert des Dirac-Operators einer kompakten Riemannschen Mannigfaltigkeit nichtnegativer Skalarkr\"{u}mmung, Math. Nachr. 97 (1980), 117-146.   
\bibitem[6]{} Th. Friedrich, On the conformal relation between twistors and Killing spinors, Supplemento de Rendiconti des Circole Mathematico de Palermo, Serie II, No.22 (1989), 59-75.
\bibitem[7]{} Th. Friedrich, Dirac-Operatoren in der Riemannschen Geometrie, Vieweg-Verlag, Braunschweig/Wiesbaden 1997. 
\bibitem[8]{} Th. Friedrich, On the spinor representation of surfaces in Euclidean 3-space, Journ. Geom. Phys. 28 (1998) 143-157.
\bibitem[9]{} O. Hijazi, A conformal lower bound for the smallest eigenvalue of the Dirac operator and Killing spinors, Comm. Math. Phys. 104, 151-162 (1986).
\bibitem[10]{} O. Hijazi, Lower bounds for eigenvalues of the Dirac operator through modified connections, Habilitation at the University of Nantes 1992; Journ. Geom. Phys. 16 (1995), 27-38.
\bibitem[11]{} E.C. Kim and Th. Friedrich, The Einstein-Dirac equation on Riemannian spin manifolds, math. DG/9905095, to appear in ''Journ. Geom. Phys.''.
\bibitem[12]{} J.M. Lee and T.H. Parker, The Yamabe problem, Bull. Amer. Math. Soc. vol.17, No.1, 1987.
\bibitem[13]{} A. Lichnerowicz, Spin manifolds, Killing spinors and the universality of the Hijazi inequality, Lett. Math. Phys. 13 (1987), 331-344.
\bibitem[14]{} J. Lott, Eigenvalue bounds for the Dirac operator, Pac. Journ. of Math. 125 (1986)
\end{Literatur}

\end{document}